\documentclass[10pt,a4paper]{article}
\usepackage{fullpage}
\usepackage{amsfonts,amsmath,amssymb}
\usepackage{amsthm}
\usepackage{graphicx}\usepackage{sectsty}

\theoremstyle{plain}
\newtheorem{theorem}{Theorem}[section]

\theoremstyle{definition}

\theoremstyle{remark}

\numberwithin{equation}{section}

\sectionfont{\large}

\begin{document}
\title{A Sub-Density Theorem of Sturm-Liouville Eigenvalue Problem with Finitely Many Singularities}
\author{Lung-Hui Chen$^1$}\maketitle\footnotetext[1]{Department of
Mathematics, National Chung Cheng University, 168 University Rd.
Min-Hsiung, Chia-Yi County 621, Taiwan. Email:
mr.lunghuichen@gmail.com. Fax:
886-5-2720497.}
\begin{abstract}
We study the distribution of the Sturm-Liouville eigenvalues of a potential with finitely many singularities. There is an asymptotically periodical structure on this class of eigenvalues as described by the entire function theory. We describe the singularities of its potential function explicitly in its eigenvalue asymptotics.
\end{abstract}
\section{Introduction and Main Result}
In this short note, we study the eigenvalue distribution for the following differential equation.
\begin{eqnarray}\label{1.1}
\left\{
  \begin{array}{ll}
    -y''(x)+p(x)y(x)=\omega^2y(x), & 0< x<\pi;\vspace{8pt}\\\vspace{8pt}
    y(0;\omega)=0,\,y'(0;\omega)=1;\\
   y(\pi;\omega)=0,
  \end{array}
\right.
 \end{eqnarray}
where
\begin{equation}
p(x)=\sum_{m=0}^Mp_m(x)+r_M(x)
\end{equation}
with $r_M(x)\in \mathcal{C}^M[0,\pi]$; Most importantly,
\begin{eqnarray}
&&p_m(x)=\sum_k c_{m,k}1_{[x_{m,k},\infty)}(x)(x-x_{m,k})^m/m!,\,m\geq1;\\
&&p_0(x)=\sum_k c_{0,k}1_{[x_{0,k},\infty)}(x),
\end{eqnarray}
where $\{x_{m,k}\}_{m,k}\in(0,\frac{\pi}{2})$ and $\{c_{m,k}\}_{m,k}\in\mathbb{R}$. We are dealing with a piecewise $\mathcal{C}^M[0,\pi]$ potential function $p(x)$. For each $m$, $x_{m,k}$ are distinct and $p(x)$ has a jump at m-th derivative at $x_{m,k}$.  We assume nontrivially the $\{x_{m,k}\}_{m,k}\in(0,\frac{\pi}{2})$ has $J$ elements and are all distinct. If there are two singular points symmetrically to the middle point located in $(0,\pi)$, then our method doesn't apply in this case.

\par
It is asked by Carlson, Threadgill and Shubin \cite{Carlson1}: How are the singularities of $p$  manifested in the distribution of eigenvalues? Being considered as a function of $\omega$, $y(\pi;\omega)$ is an entire function of $\omega$. Moreover, the zeros of $y(\pi;\omega)$ are the Dirichlet eigenvalues of the system~(\ref{1.1}). To study the asymptotics of Dirichlet eigenvalues, we examine the zeros of entire function $y(\pi;\omega)$. We try to answer the question from the point of view of complex analysis in this particular setting. In \cite{Carlson1}, a distribution of the eigenvalues  with coefficients in terms of spectral invariants is described in \cite[Theorem 4.4]{Carlson1} applying the Newton's method. In this paper, we try to characterize the distribution of the eigenvalues explicitly in terms of the singularities themselves and find the composites of the Dirichlet eigenvalues.  Can one really hear the singularities of the potential $p$? We state the main result of this paper:
\begin{theorem}\label{11}
Let $\{\omega_j\}_{j=1}^J$ be the rearrangement of the singular points $\{x_{m,k}\}_{m,k}$ such that $0=:\omega_0<\omega_{1}<\omega_2<\ldots<\omega_J<\omega_{J+1}:=\frac{\pi}{2}$. There exist exactly $2J+2$ subsequences of the zeros of $y(\pi;\omega)$, denoted as $\{z_{n_l}\}$, where $l=1,2,\ldots,2J+2$, such that
\begin{eqnarray}
&&z_{n_l}\sim\frac{{n_l}\pi}{\omega_{l}-\omega_{l-1}}+O(1),\mbox{ as }{n_l}\rightarrow\pm\infty\mbox{ in }\mathbb{Z};\\
&&\bigcup_{l=1}^{2J+2}\{z_{n_l}\}=\{z_n\},
\end{eqnarray}
in which $\{z_n\}$ are the zeros of $y(\pi;\omega)$.
\end{theorem}
In particular, we recover the point set $\{\omega_j\}_{j=1}^J$ from the subsequences of Dirichlet eigenvalues corresponding to each of these points. We may refine the asymptotics~(\ref{11}) to next order by the method in \cite[p.\,37]{Po}:
\begin{equation}
z_{n_l}\sim\frac{{n_l}\pi}{\omega_{l}-\omega_{l-1}}+O(\frac{1}{{n_l}}),\mbox{ as }{n_l}\rightarrow\pm\infty\mbox{ in }\mathbb{Z}.
\end{equation}
This is the only eigenvalue asymptotics containing the information on the position of the singularities of a given potential function known to the author. We may compare the result in \cite{Hy,Marchenko,Po}. However, in \cite{Hy}, they considered a much general class of potential functions. One may sum up all of the subsequences to obtain the classic eigenvalue density as in  \cite{Marchenko,Po}.

\par
We start with the asymptotic expansion of the solution of~(\ref{1.1}) which we refer to \cite{Carlson1,Carlson2}. The following asymptotics holds:
\begin{eqnarray}\nonumber
y(\pi;\omega)&=&\frac{\sin\{\omega\pi\}}{\omega}-\frac{1}{2}\frac{\cos\{\omega\pi\}}{\omega^2}\int_0^\pi p(t)dt\\
&&+2\sum_{m=1}^{\frac{M-1}{2}}(-1)^m[2\omega]^{-2m}\cos\{\omega\pi\}P_{m-1}(\pi)\nonumber\\\nonumber
&&+2\sum_{m=1}^{\frac{M-1}{2}}(-1)^{m+1}[2\omega]^{-2m-1}\sin\{\omega\pi\}Q_{m-1}(\pi)\\\nonumber
&&+2\sum_{m=0}^{\frac{M-1}{2}-1}(-1)^{m}[2\omega]^{-2m-3}
\times\{\sum_{x_{2m,k}<\pi}c_{2m,k}\sin\{\omega[x-2x_{2m,k}]\}\}\\\nonumber
&&+2\sum_{m=0}^{\frac{M-1}{2}-1}(-1)^{m+1}[2\omega]^{-2m-4}
\times\{\sum_{x_{2m+1,k}<\pi}c_{2m+1,k}\cos\{\omega[x-2x_{2m+1,k}]\}\}\\
&&+O(\frac{\exp{\omega\pi}}{\omega^M}),\label{1.5}
\end{eqnarray}
if $\omega\in\mathbb{C}$ and
\begin{eqnarray}
&&P_0(x):=\int_0^xp(t_1)dt_1;\\
&&P_m(x):=q_{2m-1}(p;x)-q_{2m-1}(p;0),\,m\in\mathbb{N};\\
&&Q_m(x):=q_{2m}(p;x)+q_{2m}(p;0),\,m\in\mathbb{N}_0;\\
&&q_m(p;x):=[\sum_{k=m}^Mp_k(x)+r_M(x)]^{(m)},\,r_M(x)\in\mathcal{C}^M[0,\pi].
\end{eqnarray}
This is essentially the (3.e) in \cite[p.\,84]{Carlson1}. However, we deal with $\omega\in\mathbb{C}$ in this paper. The only difference is in the big O-term in the end of~(\ref{1.5}). We refer the proof to \cite[p.\,84]{Carlson1}, and also \cite{Po}, which comes from the repeated integration by parts.  We will apply the Wilder's theorem to~(\ref{1.5}) which is a sum of asymptotically hyperbolic series to obtain the asymptoics of the Dirichlet eigenvalues.

\section{The Wilder's theorem}

There is an asymptotic periodic structure \cite{Dickson,Dickson2,Levin} within the zero set of the asymptotically hyperbolic sum, say, the asymptotic expansion~(\ref{1.5}). We refer to \cite{Dickson2,Levin} for a comprehensive study on the zero distribution theory of this kind. To be more convincing, we start with the following theorem. One can bypass this part if familiar with the entire function theory.
The indexing in this section is independent of the others.

\begin{theorem}[Dickson \cite{Dickson}] Let
\begin{equation}
R(\alpha,s,h):=\{z=x+iy\in\mathbb{C}|\,|x|\leq
h,\,y\in[\alpha,\alpha+s]\};
\end{equation}
\begin{equation}
N_g(R(\alpha,s,h)):=\{\mbox{ the number of zeros of }g(z)\mbox{ in
}R(\alpha,s,h)\},
\end{equation}
in which $$g(z)=\sum_{j=1}^nA_j e^{\omega_jz},$$ where $z=x+iy$,
$A_j\neq0$, $\omega_1<\omega_2<\cdots<\omega_n$. Then, there
exists $K>0$ such that
\begin{enumerate}
    \item each zero of $g$ is in $|x|<K$;
\item for each pair of reals $(\alpha,s)$ with $s>0$,
\begin{equation}
|N_g(R(\alpha,s,K))-s(\omega_n-\omega_1)/(2\pi)|\leq n-1.
\end{equation}
\end{enumerate}
\end{theorem}
Let us acquire a more sophisticated theorem of this type.
Let
\begin{equation}\label{2.52}
f(z)=\sum_{j=1}^nA_jz^{m_j}[1+\epsilon(z)]e^{\omega_jz},
\end{equation}
where $n>1$ and $A_j$ and $\omega_j$ are complex numbers such that
$A_j\neq0$ and the $\omega_j$ are distinct; the $m_j$ are
non-negative integers; the functions $\epsilon$ are analytic for
$|z|\geq r_0\geq0$ with $\lim_{z\rightarrow\infty}\epsilon(z)=0$.
When we are talking about the zeros of $f(z)$, we are referring to
its zeros outside certain open ball around the origin.

\par
We set up
the following quantities to the $f(z)$ in~(\ref{2.52}): Let $Q$ be
the broken line given by the $\overline{\omega}_j$ given
in~(\ref{2.52}) with
$\overline{\omega}_1,\cdots,\overline{\omega}_\sigma$ as its
vertices. The indices are labeled counterclockwise. Let $L_k$ be
the line segment $[\overline{\omega}_k,\overline{\omega}_{k+1}]$
and
$$\phi_k:=\arg\{\overline{\omega}_k-\overline{\omega}_{k+1}\}$$ in
$[-\frac{\pi}{2},\frac{3\pi}{2})$. Let
\begin{equation}
e_k=e^{i\phi_k}. \label{4.6}
\end{equation}
Certain $\overline{\omega}_p$ on $L_k$ are assigned doubly indexed
subscripts as follows: Let the convex hull of
$\overline{\omega}_{k},\,\overline{\omega}_{k+1}$ and
$\tau_p=\overline{\omega}_p+im_pe_k$ in which
$\overline{\omega}_p$ on $L_k$; assign subscripts
$j=1,\cdots,\sigma_k$ to $\omega_{kj}$ so that
$\omega_{k1}=\omega_k,\,\omega_{k\sigma_k}=\omega_{k+1}$ and
$\tau_{kj}$ are vertices of this convex hull and preceding in a
counterclockwise direction from $\overline{\omega}_k+im_ke_k$ to
$\overline{\omega}_{k+1}+im_{k+1}e_k$. For
$j=1,\cdots,\sigma_{k}-1$,
$$L_{kj}:=[\tau_{kj},\tau_{kj+1}];$$
\begin{equation} \label{2.6}
\mu_{kj}:=\frac{m_{kj}-m_{kj+1}}{(\omega_{kj}-\omega_{kj+1})e_k}
\end{equation} which is real; $n_{kj}$ is the number of $\tau_p$ on $L_{kj}$. In particular, if $L_{kj}$ in an interval with exactly two end points, then we have $n_{kj}=2$.
\par
Moreover, for  $j=1,\cdots,\sigma_{k}-1$ and $h>0$, we define
\begin{equation}\label{2.7}
V_{kj}(h):=\{z|\,\Im(z/e_k)\geq0,\,|\Re(z/e_k)+\mu_{kj}\log|z||\leq
h\}.
\end{equation}
$T_k(\theta)$ is defined to be a closed sector with vertex at zero
of opening $2\theta$ about the outward normal to $L_k$ through the
origin. For the same $k$ and $j$ and each triple of reals
$(\alpha,s,h)$, $s>0$ and $h>0$, the set
\begin{equation}\label{2.8}
R_{kj}(\alpha,s,h):=\{z|\,\Im(z/e_k)+\mu_{kj}\arg
z\in[\alpha,\alpha+s],\,|\Re(z/e_k)+\mu_{kj}\log|z||\leq h\},
\end{equation}
where $\arg z\in(\phi_k,\phi_k+\pi)$ and $R_{kj}(\alpha,s,h)$ is in $V_{kj}(h)\cap T_k(\theta)$. They are asymptotically logarithmic tubular neighborhoods. We refer to \cite{Dickson} for a comprehensive study. Now we state the
following theorem.
\begin{theorem}[Dickson \cite{Dickson}]\label{D}
Let $f(z)$ be given as in~(\ref{2.52}). Then, there exists $h>0$
such that \begin{enumerate}
            \item all but a finite number of zeros of $f$ of modulus
greater than $r_0$ are in $\bigcup_{k,j}V_{kj}$;
            \item  for each pair
of positive reals $\epsilon$ and $s_0$, there exists an
$\alpha_0=\alpha_0(\epsilon,s_0)$ such that whenever
$\alpha\geq\alpha_0$ and $s\geq s_0$,
\begin{equation}\label{2.30}
|N_f(R_{kj}(\alpha,s,h))-s|\omega_{kj+1}-\omega_{kj}|/(2\pi)|<n_{kj}-1+\epsilon.
\end{equation}
\end{enumerate}
\end{theorem}
This is exactly stated as in \cite{Dickson}. The proof is
in \cite[Theorem 2, p.21]{Dickson2}. We refer to \cite{Chen} for another application of this theorem.

\section{Proof of Theorem \ref{11}}
\begin{proof}
To apply Theorem \ref{D} to the hyperbolic sum~(\ref{1.5}), we rewrite~(\ref{1.5}): It is well-known \cite{Po} that
there is a $C_\delta$ depending on the distance to the zeros of $\sin{\omega\pi}$ such that
\begin{equation}
\exp|\Im\omega\pi|<C_\delta\sin\{\omega\pi\}.
\end{equation}
Hence,~(\ref{1.5}) becomes
\begin{eqnarray}\nonumber
y(\pi;\omega)&=&\frac{\sin\{\omega\pi\}}{\omega}[1+O(\frac{C_\delta}{\omega^{M-1}})]
-\frac{1}{2}\frac{\cos\{\omega\pi\}}{\omega^2}\int_0^\pi p(t)dt\\
&&+2\sum_{m=1}^{\frac{M-1}{2}}(-1)^m[2\omega]^{-2m}\cos\{\omega\pi\}P_{m-1}(\pi)\nonumber\\\nonumber
&&+2\sum_{m=1}^{\frac{M-1}{2}}(-1)^{m+1}[2\omega]^{-2m-1}\sin\{\omega\pi\}Q_{m-1}(\pi)\\\nonumber
&&+2\sum_{m=0}^{\frac{M-1}{2}-1}(-1)^{m}[2\omega]^{-2m-3}
\times\{\sum_{x_{2m,k}<\pi}c_{2m,k}\sin\{\omega[x-2x_{2m,k}]\}\}\\\nonumber
&&+2\sum_{m=0}^{\frac{M-1}{2}-1}(-1)^{m+1}[2\omega]^{-2m-4}
\times\{\sum_{x_{2m+1,k}<\pi}c_{2m+1,k}\cos\{\omega[x-2x_{2m+1,k}]\}\},\,\omega\notin\mathbb{Z}.
\end{eqnarray}
Now we rearrange according to their exponential powers by the theorem assumption to the following form:
\begin{eqnarray}\nonumber
y(\pi;\omega)&=&\{\frac{-1}{2i\omega}[1+O(\frac{C_\delta}{\omega^{M-1}})]
+\sum_{m=1}^{\frac{M-1}{2}}(-1)^m[2\omega]^{-2m}P_{m-1}(\pi)\\\nonumber
&&+i\sum_{m=1}^{\frac{M-1}{2}}(-1)^{m+1}[2\omega]^{-2m-1} Q_{m-1}(\pi)\}e^{-i\omega\pi}\\\nonumber
&&+\sum_{j=1}^J C_j(m,\omega)e^{-i\omega(\pi-2\omega_j)}+\sum_{j=J}^1 D_j(m,\omega)e^{i\omega(\pi-2\omega_j)}\\
&&+\{\frac{1}{2i\omega}[1+O(\frac{C_\delta}{\omega^{M-1}})]+\sum_{m=1}^{\frac{M-1}{2}}(-1)^m[2\omega]^{-2m}P_{m-1}(\pi)
\nonumber\\&&-i\sum_{m=1}^{\frac{M-1}{2}}(-1)^{m+1}[2\omega]^{-2m-1} Q_{m-1}(\pi)\}e^{i\omega\pi},\,\omega\notin\mathbb{Z}.\label{3.1}
\end{eqnarray}
in which the $C_j(m,\omega)$ and $D_j(m,\omega)$ can be obtained by comparing~(\ref{3.1}) with~(\ref{1.5}). Besides~(\ref{3.1}), the entire function $y(\pi;\omega)$ is bounded near $\mathbb{Z}$.
Without loss of generality, we consider the zeros of $$Y(\omega):=\omega y(\pi;\omega/i)$$ by applying Theorem \ref{D} in a suitable strip containing the real axis. We observe the zeros of $Y(\omega)$ spread themselves vertically along the imaginary axis,  that is , the zeros of $y(\pi;\omega)$ spread themselves along the real axis. In particular,  given the singularity sequence $\{x_{m,k}\}_{m,k}$ which are all distinct by assumption with $J$ elements, we let $\{\omega_j\}_{j=1}^J$ be the rearrangement of $\{x_{m,k}\}_{m,k}$ such that $$0=\omega_0<\omega_{1}<\omega_2<\ldots<\omega_J<\omega_{J+1}=\frac{\pi}{2}.$$ We construct following $2J+2$ successive intervals in $[-\pi,\pi]$: \begin{eqnarray*}
&&L_1:=[-\pi,-\pi+2\omega_1],\,L_2:=[-\pi+2\omega_1,-\pi+2\omega_2],\ldots,\,L_{J+1}:=[-\pi+2\omega_J,0],\\
&&L_{J+2}:=[0,\pi-2\omega_J],\ldots,
L_{2J+1}:=[\pi-2\omega_2,\pi-2\omega_1],\,L_{2J+2}:=[\pi-2\omega_1,\pi].
\end{eqnarray*}
These intervals are applied as the polygons described previously. However, we note that $L_{J+1}\cup L_{J+2}$ combines to generate a sequence of zeros as described by~(\ref{2.8}) and then~(\ref{2.30}) after observing the exponential exponents in~(\ref{3.1}). There are actually $2J+1$ asymptotically rectangular area on duty. Without loss of generality, we take each $\{L_{l}\}_{l=1}^{2J+2}$ to generate an asymptotically rectangular area $\{R_l\}_{l=1}^{2J+2}$ as described by~(\ref{2.7}) and~(\ref{2.8}). Finally, we note that one can not identify the quantities $\{\mu_{k,j}\}$ in~(\ref{2.6}), because not being able to locate the coefficients $\{m_j\}$ in~(\ref{2.52}) again in~(\ref{3.1}). For our case, the quantities $\{e_k\}$ in~(\ref{4.6}) are equal to $1$.

\par
Let zeros in $R_l$ be denoted as $z_{n_l}$, $l=1,2,\ldots,2J+2$ and $n_{l}\in\mathbb{N}$.
Hence,~(\ref{2.30}) implies that
\begin{equation}\label{3.2}
|N_Y(R_l(\alpha,s,h))-s(\omega_{l}-\omega_{l-1})/\pi|<\tilde{n}_l-1+\epsilon,\,l=1,\ldots,2J+2,
\end{equation}
in which $\tilde{n}_l=2$ by the construction of intervals $\{L_l\}_{l=1}^{2J+2}$ for all $l$. Here, we define $\omega_l$, $J+1\leq l\leq2J+2$ by the symmetry of $\{L_l\}_{l=1}^{2J+2}$. Finally, it is well-known that the zeros of $y(\pi;\omega)$ are simple and real \cite{Po}, so~(\ref{3.2}) implies
\begin{equation}\nonumber
z_{n_l}\sim\frac{{n_l}\pi}{\omega_{l}-\omega_{l-1}}+O(1),\,l=1,\ldots,2J+2;\,n_{l}\in\mathbb{N}.
\end{equation}
Because the zeros of $y(\pi;\omega)$ are symmetric to the imaginary axis, we can rewrite the equation above to be
\begin{equation}\label{34}
z_{n_l}\sim\frac{{n_l}\pi}{\omega_{l}-\omega_{l-1}}+O(1),\,l=1,\ldots,2J+2;\,n_{l}\in\mathbb{Z}.
\end{equation}
Once again, we note that $\{z_{n_{J+1}}\}\cup\{z_{n_{J+2}}\}$ is the sequence of zeros generated by the interval $L_{J+1}\cup L_{J+2}$.
This proves the Theorem \ref{11}.
\end{proof}
We may observe from~(\ref{3.2}) that the total number of the zeros of $y(\pi;\omega)$ is equal to
\begin{eqnarray}
N_Y(\bigcup_{l=1}^{2J+2}R_{l}(\alpha,s,h))
=\sum_{l=1}^{2J+2}N_Y(R_{l}(\alpha,s,h))
\sim \sum_{l=1}^{2J+2}s(\omega_{l}-\omega_{l-1})/\pi+o(1)
=s+o(1),
\end{eqnarray}
which implies the zeros of $y(\pi;\omega)$ denoted as , $\{z_n\}_{n\in\mathbb{Z}}$, have the following asymptotics:
\begin{equation}
z_n\sim n+O(1),\,n\in\mathbb{Z},
\end{equation}
which matches with the Counting Lemma in \cite[p.\,36]{Po}. The bounded error term in~(\ref{34}) avoids the possible contradiction to the classic asymptotics \cite{Po}.

\end{document}